\documentclass[12pt]{article}
\usepackage{dsfont}
\usepackage{latexsym,amsfonts,amssymb}
\usepackage[dvips]{color}
\usepackage{colordvi,multicol}

\newtheorem{thm}{Theorem}[section]
\newtheorem{cor}[thm]{Corollary}

\newtheorem{defi}{Definition}[section]
\newtheorem{rem}{Remark}[section]

\begin{document}
\title{\bf  A Note on Uniform Nonintegrability of Random Variables}
\author{Ze-Chun Hu and Xue Peng\thanks{Corresponding author.\ E-mail addresses: zchu@scu.edu.cn (Z.-C. Hu), pengxuemath@scu.edu.cn (X. Peng).}\\
 {\small College of Mathematics, Sichuan  University, Chengdu, China}}
\date{}
\maketitle

\begin{abstract}

In  a recent paper \cite{CHR16}, Chandra, Hu and Rosalsky introduced the notion of a sequence of random variables being uniformly nonintegrable, and presented  a list of interesting results on this uniform nonintegrability. In this note,  we introduce a weaker definition on uniform nonintegrability (W-UNI for short) of random variables,  present a
necessary and sufficient condition for W-UNI, and  give two equivalent characterizations of W-UNI, one
of which is a W-UNI analogue of the celebrated de La Vall\'{e}e Poussin criterion
for uniform integrability. In addition, we give some remarks, one of which gives a negative answer to the open problem raised  in \cite{CHR16}.

\end{abstract}



\noindent

\section{Introduction}

It's well known that the uniform integrability of a sequence of random variables plays an important role in probability theory. As to the uniform integrability criterions, please refer to  Chung (1974, P. 96), Chong (1979),  Chow and Teicher (1997, P. 94), Hu and Rosalsky (2011), Klenke (2014, p. 138) and Chandra (2015).

In a recent paper \cite{CHR16}, the authors introduced the notion of a sequence of random variables being {\it
uniformly nonintegrable} and gave some interesting characterizations of this uniform nonintegrability. Now, we recall that definition in \cite{CHR16}. Let $(\Omega,\mathcal{F},P)$ be a probability space. Suppose that all random variables under consideration are defined on this probability space. Let $X$ be a random variable and $A\in \mathcal{F}$. We denote $E(XI_A)$ by $E(X: A)$.

For a random variable $X$, it's said to be {\it nonintegralbe} (NI) if $E|X|=\infty$. By the Lebesgure monotone convergence theorem, we have
\begin{eqnarray}\label{1}
X\ \mbox{is NI if and only if}\ \lim_{a\to \infty}E(|X|: |X|\leq a)=\infty.
\end{eqnarray}

In virture of (\ref{1}),  the authors of \cite{CHR16} gave the following definition:

\begin{defi}\label{defi-1.1}
A sequence of random variables $\{X_n,n\geq \}$ is said to be {\it uniformly nonintegrable} (UNI) if
\begin{eqnarray}\label{defi1.1-a}
\lim_{a\to\infty}\inf_{n\geq 1}E(|X_n|: |X_n|\leq a)=\infty
\end{eqnarray}
or equivalently,
$$
\sup_{N\geq 1}\inf_{n\geq 1}E(|X_n|: |X_n|\leq N)=\infty.
$$
\end{defi}

By the Lebesgure monotone convergence theorem, we also have
\begin{eqnarray}\label{2}
X\ \mbox{is NI if and only if}\ \lim_{a\to\infty}E(|X|\wedge a)=\infty.
\end{eqnarray}

In virture of (\ref{2}), we give the following definition.

\begin{defi}\label{defi-1.2}
A sequence of random variables $\{X_n,n\geq \}$ is said to be {\it W-uniformly nonintegrable} (W-UNI) if
\begin{eqnarray}\label{defi1.1-a}
\lim_{a\to\infty}\inf_{n\geq 1}E(|X_n|\wedge a)=\infty
\end{eqnarray}
or equivalently,
$$
\sup_{N\geq 1}\inf_{n\geq 1}E(|X_n|\wedge N)=\infty,
$$
where ``W" stands for ``weak".
\end{defi}

In fact, for any random varialbe $X$ and  positive constant $a$, we have
$
E(|X|\wedge a)\geq E(|X|: |X|\leq a).
$
Hence UNI $\Rightarrow$ W-UNI.

\begin{rem}\label{rem1.1-a}
If $\{X_n,n\geq \}$ is W-UNI, then $X_n$ is NI for all $n\geq 1$. In fact, let $k\geq 1$. Then for $a>0$,
$$
E|X_k|\geq E(|X_k|\wedge a)\geq \inf_{n\geq 1}E(|X_n|\wedge a)\to \infty\ \mbox{as}\ a \to \infty.
$$
\end{rem}

\begin{rem}\label{rem1.1}
Let $\{X_n,n\geq 1\}$ be a sequence of random variables such that $|X_n|\geq Y$ a.s., $n\geq 1$, where $Y$ is a NI random variable. By Definition \ref{defi-1.2}, we can easily get that $\{X_n,n\geq 1\}$ is W-UNI. But \cite[Example 4.1]{CHR16} shows that  we cann't conclude that $\{X_n,n\geq 1\}$ is UNI.
\end{rem}

In Section 2, we will give two examples to show that UNI is strictly stronger than W-UNI in general, and even the condition (3.1) defined in \cite[Theorem 3.1]{CHR16} fails, W-UNI can still holds. In Section 3, we present a
necessary and sufficient condition for W-UNI, and  give two equivalent characterizations of W-UNI, one
of which is a W-UNI analogue of the celebrated de La Vall\'{e}e Poussin criterion
for uniform integrability. In the final section, we give some remarks, one of which gives a negative answer to the open problem raised in \cite[Remark 5.3]{CHR16}.

\section{Two examples}\setcounter{equation}{0}

In this section, we borrow two examples from \cite{CHR16} to illustrate some phenomena. The first example shows that  UNI is strictly stronger than W-UNI in general.

\noindent {\bf Example 2.1} (\cite[Example 4.1]{CHR16}) Let $Y$ be a random variable with probability density function
$$
f(y)=\frac{1}{y^2}I_{[1,\infty)}(y),\ -\infty<y<\infty
$$
and set $X_n=nY,n\geq 1$. The authors of \cite{CHR16} showed that $\{X_n,n\geq 1\}$ is not UNI.

Obviously, we have that for any $n\geq 1, X_n\geq Y\geq 0$ a.s., and $Y$ is a NI random variable.  By Remark \ref{rem1.1}, we get that $\{X_n,n\geq 1\}$ is  M-UNI, which can be verified directly as follows. For any $a>0$,  we have by direct calculation that
\begin{eqnarray*}
E[|X_n|\wedge a]=\left\{
\begin{array}{cl}
a,& \mbox{if}\ n\geq a;\\
n(1+\ln \frac{a}{n}),& \mbox{if}\ n<a.
\end{array}
\right.
\end{eqnarray*}
For any $a>1$, define a function
$$
f(x)=x\left(1+\ln \frac{a}{x}\right),\ 1\leq x<a.
$$
Then we have
$$
f'(x)=\ln a-\ln x>0,\ \forall x\in [1,a),
$$
and thus $\inf_{1\leq x<a}f(x)=f(1)=1+\ln a$. Hence we have
\begin{eqnarray*}
\lim_{a\to\infty}\inf_{n\geq 1}E[|X_n|\wedge a]\geq \lim_{a\to\infty}(a\wedge (1+\ln a))=\infty,
\end{eqnarray*}
i.e. $\{X_n,n\geq 1\}$ is  W-UNI.

\bigskip

\cite[Theorem 3.1]{CHR16}(ii) says that if $\{X_n,n\geq 1\}$ is  UNI, then the condition (3.1) defined in \cite[Theorem 3.1]{CHR16} must  hold, i.e.
\begin{eqnarray}\label{2.1}
\beta:=\lim_{a\to\infty}\inf_{n\geq 1}P(|X_n|\leq a)>0.
\end{eqnarray}
The following example shows that even the above condition fails, W-UNI can still holds.

\noindent {\bf Example 2.2} (\cite[Example 4.2]{CHR16}) Let $\{a_n,n\geq 1\}$ be a sequence in $(0,1)$, and let $\{X_n,n\geq 1\}$ be a sequence of random variables where for each $n\geq 1, X_n$ has probability density function
$$
f_n(x)=\frac{1-\alpha_n}{x^{2-\alpha_n}}I_{[1,\infty)}(x),\ -\infty<x<\infty.
$$
We suppose that $\sup_{n\geq 1}\alpha_n=1$ and $\inf_{n\geq 1}\alpha_n=0$. Then by \cite[Example 4.2]{CHR16}, we know that
$$
\beta=\lim_{n\to\infty}\inf_{n\geq 1}\left(1-\frac{1}{a^{1-\alpha_n}}\right)=\lim_{n\to\infty}\left(1-\frac{1}{a^{1-\sup_{n\geq 1}\alpha_n}}\right)=0.
$$

For $a>1$ and $n\geq 1$, by \cite[Example 4.2]{CHR16}, we have
\begin{eqnarray*}
E[|X_n|\wedge a]&=&E[|X_n|: |X_n|\leq  a]+aP(|X_n|>a)\\
&=&\frac{1-\alpha_n}{\alpha_n}(a^{\alpha_n}-1)+a\int_a^{\infty}\frac{1-\alpha_n}{x^{2-\alpha_n}}dx\\
&=&\frac{1-\alpha_n}{\alpha_n}(a^{\alpha_n}-1)+a^{\alpha_n}\\
&=&\frac{1}{\alpha_n}(a^{\alpha_n}-1)+1.
\end{eqnarray*}

For $a>3$, define a function
\begin{eqnarray*}
g(x)=\left\{
\begin{array}{cl}
\frac{1}{x}\left(a^x-1\right),& \mbox{if}\quad 0<x\leq 1;\\
\ln a,& \mbox{if}\quad x=0.
\end{array}
\right.
\end{eqnarray*}
By L'Hospital principle, we have
$$\lim_{x\to 0+}g(x)=\lim_{x\to 0+}a^x\ln a=\ln a.$$
So  $g(x)$ is a continuous function on $[0,1]$.

For any $0<x<1$, we have
$$
g'(x)=\frac{1}{x^2}+\frac{a^x}{x}\left(\ln a-\frac{1}{x}\right)=\frac{1}{x}\left[a^x\ln a-\frac{a^x-1}{x}\right].
$$
By the mean value theorem, we obtain that there exists $\theta\in (0,x)$ such that
$\frac{a^x-1}{x}=a^{\theta}\ln a$. So for any $x\in (0,1)$, we have
$$
g'(x)=\frac{1}{x}\left[(a^x-a^{\theta})\ln a\right]>0.
$$
Hence $g(x)$ is an strictly increasing function on $[0,1]$. By the assumption  $\inf_{n\geq 1}\alpha_n=0$ and L'Hospital principle, we have
\begin{eqnarray*}
\lim_{a\to\infty}\inf_{n\geq 1}E[|X_n|\wedge a]&=& \lim_{a\to\infty}\inf_{n\geq 1}\left(\frac{1}{\alpha_n}(a^{\alpha_n}-1)+1\right)\\
&=&\lim_{a\to\infty}\lim_{\alpha\downarrow 0}\left(\frac{1}{\alpha}(a^{\alpha}-1)+1\right)\\
&=&\lim_{a\to\infty}(\ln a+1)=\infty,
\end{eqnarray*}
i.e. $\{X_n,n\geq 1\}$ is  W-UNI.

\section{The main results}\setcounter{equation}{0}

At first, we  give a necessary and sufficient condition for a sequence of random variables to be  W-UNI.

\begin{thm}\label{thm-equivalent condition for M-UNI} Let $\{X_n, n\geq1\}$ be a sequence of random variables. Then $\{X_n, n\geq1\}$ is  M-UNI if and only if for all $M>0$, there exists $\alpha\in (0,1)$ such that for every sequence of events $\{A_n,n\geq 1\}$,
$$
\inf_{n\geq 1}P(A_n)\geq \alpha \Rightarrow \inf_{n\geq 1}E(|X_n|: A_n)\geq M.
$$
\end{thm}
\noindent {\bf Proof.} (Necessity) Suppose that $\{X_n, n\geq1\}$ is M-UNI, i.e. $\lim\limits_{a\to\infty}\inf_{n\geq1}E(|X_n|\wedge a)=\infty$. Then for all $M>0$, there exsits $a_0>M$ such that
$\inf_{n\geq1}E(|X_n|\wedge a_0)>4M$.
Hence, for any $n\geq 1$,
\begin{eqnarray*}
E(|X_n|:|X_n|\leq a_0)+a_0P(|X_n|>a_0)=E(|X_n|\wedge a_0)>4M,
\end{eqnarray*}
which implies either $E(|X_n|:|X_n|\leq a_0)>2M$ or $a_0P(|X_n|>a_0)>2M$.  Let $\alpha=1-\frac{M}{a_0}$, and a sequence of events $\{A_n,n\geq1\}$ satisfying $\inf_{n\geq1}P(A_n)\geq\alpha$.

When $E(|X_n|:|X_n|\leq a_0)>2M$,  we have
\begin{eqnarray*}
E(|X_n|:A_n)
&\geq&E(|X_n|:[|X_n|\leq a_0]\cap A_n)\\
&=&E(|X_n|:|X_n|\leq a_0)-E(|X_n|:[|X_n|\leq a_0]\cap A^c_n)\\
&>&2M-a_0P(A^c_n)\\
&=&2M+a_0(P(A_n)-1)\\
&\geq&2M+a_0(\alpha-1)=M.
\end{eqnarray*}

And when $a_0P(|X_n|>a_0)>2M$, we get
\begin{eqnarray*}
E(|X_n|:A_n)
&\geq& E(|X_n|:[|X_n|>a_0]\cap A_n)\\
&\geq& a_0P([|X_n|>a_0]\cap A_n)\\
&=&a_0[P(|X_n|>a_0)-P([|X_n|>a_0]\cap A^c_n)]\\
&\geq& a_0P(|X_n|>a_0)-a_0P(A^c_n)\\
&=& a_0P(|X_n|>a_0)+a_0(P(A_n)-1)\\
&>&2M+a_0(\alpha-1)=M.
\end{eqnarray*}
Hence $\inf_{n\geq1}E(|X_n|;A_n)\geq M$.

(Sufficiency) Suppose that  for all $M>0$, there exists $\alpha\in (0,1)$ such that for every sequence of events $\{A_n,n\geq 1\}$,
\begin{eqnarray}\label{proof-of-thm3.1}
\inf_{n\geq 1}P(A_n)\geq \alpha \Rightarrow \inf_{n\geq 1}E(|X_n|: A_n)\geq M.
\end{eqnarray}

 Note that $\inf_{n\geq1}P(|X_n|\leq a)$  is an increasing function of $a$. Then we have
\begin{center}
$\beta:=\lim\limits_{a\to\infty}\inf_{n\geq1}P(|X_n|\leq a)\geq 0$.
\end{center}

Now we have two cases:

(1) $\alpha\in(0,\beta)$. In this case, $\beta>0$. Then by the definition of $\beta$,  there exists $a_0$ such that for all $a\geq a_0$, $\inf_{n\geq1}P(|X_n|\leq a)\geq\alpha$. Hence by (\ref{proof-of-thm3.1}), we have
\begin{center}
$\inf_{n\geq1}E(|X_n|\wedge a)\geq\inf_{n\geq1}E(|X_n|:|X_n|\leq a)\geq M.$
\end{center}

(2)  $\alpha\in[\beta,\infty)\cap(0,1)$  (in this case $\beta$ may be zero). Since  $\inf_{n\geq1}P(|X_n|\leq a)$ is an increasing function of $a$, we have for all $a>0$,
\begin{center}
$\inf_{n\geq1}P(|X_n|\leq a)\leq\beta\leq\alpha.$
\end{center}
 Let $a>\frac{M}{1-\alpha}$. Decompose the positive integers set $\mathds{N}$ into two subsets $\mathds{N}_1$ and $\mathds{N}_2$ such that $\mathds{N}_1\cap\mathds{N}_2=\emptyset$ and for all $n_k\in\mathds{N}_1$,$P(|X_{n_k}|\leq a)\leq\alpha$ and for all $m_j\in\mathds{N}_2$, $P(|X_{m_j}|\leq a)>\alpha$.

(i) For any $n_k\in\mathds{N}_1$,  we have $P(|X_{n_k}|>a)\geq1-\alpha$, and thus
\begin{center}
$E(|X_{n_k}|\wedge a)\geq aP(|X_{n_k}|>a)>\frac{M}{1-\alpha}(1-\alpha)=M.$
\end{center}

(ii) For any $m_j\in\mathds{N}_2$, let $A_{m_j}=\{|X_{m_j}|\leq a\}$. Then $P(A_{m_j})>\alpha.$  Set $A_{n_k}=A_{m_1}$ for all $n_k\in\mathds{N}_1$. Then
\begin{center}
$\inf_{n\geq1}P(A_n)=\inf_{j\geq1}P(A_{m_j})\geq\alpha$,
\end{center}
which together with the assumption (\ref{proof-of-thm3.1}) implies that
$$
\inf_{j\geq1}E(|X_{m_j}|:A_{m_j})\geq\inf_{n\geq1}E(|X_n|:A_n)\geq M.
$$
It follows that
$$
E(|X_{m_j}|\wedge a)\geq E(|X_{m_j}|:|X_{m_j}|\leq a)\geq M.
$$
Hence, by (i) and (ii) we know that when $a>\frac{M}{1-\alpha}$, for all $n\in\mathds{N}$,  $E(|X_n|\wedge a)\geq M,$ and thus $\inf_{n\geq1}E(|X_n|\wedge a)\geq M$.

Hence, by (1) and (2), we obtain that when $a>\frac{M}{1-\alpha}\vee a_0$, it holds that
\begin{center}
$\inf_{n\geq1}E(|X_n|\wedge a)\geq M$,
\end{center}
which means that $\lim\limits_{a\to\infty}\inf_{n\geq1}E(|X_n|\wedge a)=\infty,$ i.e. $\{X_n,n\geq 1\}$ is M-UNI.
\hfill\fbox

\begin{rem}
Let $\{X_n,n\geq 1\}$ be a sequence of random varialbes. Define the following three conditions:

(a) $\beta:=\lim_{a\to\infty}\inf_{n\geq 1}P(|X_n|\leq a)>0$.

(b) For all $M>0$, there exists $\alpha\in (0,\beta)$ such that for every sequence of events $\{A_n,n\geq 1\}$,
$$
\inf_{n\geq 1}P(A_n)\geq \alpha \Rightarrow \inf_{n\geq 1}E(|X_n|: A_n)\geq M.
$$

(c) For all $M>0$, there exists $\alpha\in (0,1)$ such that for every sequence of events $\{A_n,n\geq 1\}$,
$$
\inf_{n\geq 1}P(A_n)\geq \alpha \Rightarrow \inf_{n\geq 1}E(|X_n|: A_n)\geq M.
$$
Then \cite[Theorem 3.1]{CHR16} can be expressed by
$$
(a)+(b)\Rightarrow \{X_n,n\geq 1\} \ \mbox{is UNI}\  \Rightarrow (a)+(c),
$$
and the above Theorem \ref{thm-equivalent condition for M-UNI} can be expressed by
$$
\{X_n,n\geq 1\} \ \mbox{is W-UNI}\  \Leftrightarrow (c),
$$
which shows that the condtion (a) is irrelevant to W-UNI. 
\end{rem}

Next, we give a characterization on W-UNI (see Theorem \ref{thm3.2} below), which  corresponds to \cite[Theorem 3.2]{CHR16}. It's well known that  for any random variable $Y$,
\begin{eqnarray*}
\sum_{n=1}^{\infty}P(|Y|>n)\leq E|Y|\leq \sum_{n=0}^{\infty}P(|Y|>n).
\end{eqnarray*}
It follows that
\begin{eqnarray}\label{3}
E|Y|=\infty&\Leftrightarrow& \sum_{n=0}^{\infty}P(|Y|>n)=\infty\Leftrightarrow \lim_{m\to\infty}\sum_{n=0}^mP(|Y|>n)=\infty.
\end{eqnarray}

In virture of (\ref{3}), we  give the following definition.

\begin{defi}\label{defi-1.3}
A sequence of random variables $\{X_n,n\geq \}$ is said to be {\it W*-uniformly nonintegrable} (W*-UNI) if
\begin{eqnarray*}\label{defi1.1-a}
\lim_{m\to\infty}\inf_{k\geq 1}\sum_{n=0}^mP(|X_k|>n)=\infty
\end{eqnarray*}
or equivalently,
$$
\sup_{m\geq 1}\inf_{k\geq 1}\sum_{n=0}^mP(|X_k|>n)=\infty.
$$
\end{defi}

\begin{rem}
By
$\sum^{m}_{n=1}P(|X_k|>n)\leq\sum^{m}_{n=0}P(|X_k|>n)\leq 1+\sum^{m}_{n=1}P(|X_k|>n),$
we get
\begin{center}
$\inf_{k\geq1}\sum^{m}_{n=1}P(|X_k|>n)\leq\inf_{k\geq1}\sum^{m}_{n=0}P(|X_k|>n)\leq 1+\inf_{k\geq1}\sum^{m}_{n=1}P(|X_k|>n).$
\end{center}
Thus,
\begin{eqnarray}\label{eq: equivalent def of WUNI}
\lim_{m\to\infty}\inf_{k\geq 1}\sum_{n=0}^mP(|X_k|>n)=\infty\Longleftrightarrow \lim_{m\to\infty}\inf_{k\geq 1}\sum_{n=1}^mP(|X_k|>n)=\infty.
\end{eqnarray}
\end{rem}

\begin{thm}\label{thm3.2}
  M-UNI $\Longleftrightarrow$ W*-UNI, i.e. a sequence of random variables $\{X_n,n\geq 1\}$ is W-UNI if and only if
  $$
  \lim_{m\to\infty}\inf_{k\geq 1}\sum_{n=0}^mP(|X_k|>n)=\infty.
  $$
\end{thm}

\noindent {\bf Proof.} ``$\Rightarrow$" For any random varialbe $X$ and any positive integer $m$, by Fubini's theorem, we have
\begin{eqnarray*}
\sum_{n=0}^mP(|X|>n)&\geq& \sum_{n=0}^m\int_n^{n+1}P(|X|>x)dx\\
&=&\int_0^{m+1}P(|X|>x)dx\\
&=&\int_0^{m+1}\left(\int_{\Omega}I_{\{|X|>x\}}dP\right)dx\\
&=&\int_{\Omega}\left(\int_0^{m+1}I_{\{|X|>x\}}dx\right)dP\\
&=&E(|X|\wedge (m+1)).
\end{eqnarray*}
It follows that W-UNI $\Rightarrow$ W*-UNI.

``$\Leftarrow $" For any random varialbe $X$ and any positive integer $m$, by Fubini's theorem, we have
\begin{eqnarray*}
\sum^{m}_{n=1}P(|X|>n)
&=&\sum^{m}_{n=1}\int^{n+1}_nP(|X|>n)dx\\
&\leq&\sum^{m}_{n=1}\int^{n+1}_nP(|X|>x-1)dx\\
&=&\int^{m+1}_1P(|X|>x-1)dx\\
&=&\int^m_0P(|X|>x)dx\\
&=&\int_{\Omega}\left(\int^m_0I_{\{|X|>x\}}dx\right)dP\\
&=&E(|X|\wedge m).
\end{eqnarray*}
Hence, by  (\ref{eq: equivalent def of WUNI}), we know that   when $\lim_{m\to\infty}\inf_{k\geq 1}\sum_{n=0}^mP(|X_k|>n)=\infty$,
$\lim_{m\to\infty}\inf_{k\geq 1}E(|X_k|\wedge m)=\infty$. Hence W*-UNI $\Rightarrow$ W-UNI.
\hfill\fbox

Next, we give a characterization on W-UNI, which  corresponds to \cite[Theorem 3.3]{CHR16}, and can be regarded as a W-UNI analogue of the celebrated de La Vall\'{e}e Poussin criterion for uniform integrability.

\begin{thm}\label{thm3.3}
A sequence of random variables $\{X_m,m\geq 1\}$ is W-UNI if and only if there exists a continuous function $\varphi:[0,\infty)\to [0,\infty)$ such that $\varphi(0)=0, \varphi$ is strictly increasing, $\varphi(x)\to \infty$ as $x\to \infty,x^{-1}\varphi(x)$ is strictly decreasing to 0 as $0<x\uparrow \infty$, and such that $\{\varphi(|X_m|),m\geq 1\}$ is W-UNI.
\end{thm}

\noindent {\bf Proof.} The main idea comes from the proof of \cite[Theorem 3.3]{CHR16}. For the reader's convenience, we spell out the details.

(Sufficiency) Suppose that there exists a function $\varphi$ satisfying the advertised properties. Then, there exists $A>0$ such that for all $x>A$, $\varphi(x)\leq x$. Thus we get
\begin{eqnarray*}
\infty
&=& \lim\limits_{a\to\infty}\inf_{m\geq1}E(\varphi(|X_m|)\wedge a)\\
&=& \lim\limits_{a\to\infty}\inf_{m\geq1}[E(\varphi(|X_m|)\wedge a:|X_m|\leq A)+E(\varphi(|X_m|)\wedge a:|X_m|> A)]\\
&\leq& \lim\limits_{a\to\infty}\inf_{m\geq1}[\varphi(A)+E(|X_m|\wedge a)]\\
&=& \varphi(A)+\lim\limits_{a\to\infty}\inf_{m\geq1}E(|X_m|\wedge a).
\end{eqnarray*}
It follows that  $\{X_m,m\geq 1\}$ is W-UNI.

(Necessity) Assume that $\{X_m,m\geq 1\}$ is W-UNI.  Firstly, we show that there exists a sequence of positive integers $\{n_k,k\geq 0\}$ such that $n_{k+1}>2n_k$ for all $k\geq 0$ and for all $j\geq1$,
\begin{eqnarray}\label{eq:proof phi(x)01}
\inf_{m\geq1}\sum\limits^{n_j}_{n=n_{j-1}+1}P(|X_m|>jn)>1.
\end{eqnarray}
Let $n_0=0$.  We can easily check  that  $\{X_m,m\geq 1\}$ is W-UNI if and only if $\{\frac{X_m}{j},m\geq 1\}$ is W-UNI for all $j\geq 1$ (or for some $j\geq1$).
Then, by Theorem \ref{thm3.2} and (\ref{eq: equivalent def of WUNI}), it holds for all $j\geq1$,
\begin{eqnarray}\label{eq:proof phi(x)02}
\lim\limits_{N\to\infty}\inf_{m\geq1}\sum\limits^{N}_{n=1}P(|X_m|>jn)=\infty.
\end{eqnarray}
By (\ref{eq:proof phi(x)02}) with $j=1$, there exists $n_1>2$ such that (\ref{eq:proof phi(x)01}) holds with $j=1$.
Suppose $n_0,n_1,\dots,n_k$ have been chosen for some $k\geq 1$ such that $n_j>2n_{j-1}$ for $j=0,1,\ldots,k-1$ and (\ref{eq:proof phi(x)01}) holds for $j=1,2,\ldots, k$. For $j=k+1$, by (\ref{eq:proof phi(x)02}) and
\begin{eqnarray*}
&&\inf_{m\geq1}\sum\limits^N_{n=1}P(|X_m|>(k+1)n)\\
&=&\inf_{m\geq1}\left[\sum\limits^{n_k}_{n=1}P(|X_m|>(k+1)n)+\sum\limits^{N}_{n=n_k+1}P(|X_m|>(k+1)n)\right]\\
&\leq& n_k+\inf_{m\geq1}\sum\limits^{N}_{n=n_k+1}P(|X_m|>(k+1)n),
\end{eqnarray*}
we get
$$\lim\limits_{N\to\infty}\inf_{m\geq1}\sum\limits^{N}_{n=n_k+1}P(|X_m|>(k+1)n)=\infty.$$
Thus there exists $n_{k+1}>2n_k$ such that (\ref{eq:proof phi(x)01}) holds with $j=k+1$  thereby establishing the above assertion.

Next, define a function $g$ on $[0,\infty)$ as follows:
$$
g(x)=k+\frac{1}{n_{k+1}-n_k}(x-n_k),\,\,n_k\leq x<n_{k+1},\,\,k\geq 0.
$$
Then $g(0)=0$ and $g$ is continuous and strictly increasing with $g(x)\to\infty$ as $x\to\infty$.
Hence, by (\ref{eq:proof phi(x)01}), for all $k\geq 1, m\geq1$, we have
\begin{eqnarray}\label{3.6}
\sum\limits^{n_k}_{n=1}P(|X_m|>ng(n))
&=& \sum\limits^k_{j=1}\sum^{n_j}_{n=n_{j-1}+1}P(|X_m|>ng(n))\nonumber\\
&\geq& \sum\limits^k_{j=1}\sum^{n_j}_{n=n_{j-1}+1}P(|X_m|>ng(n_j))\nonumber\\
&=& \sum\limits^k_{j=1}\sum^{n_j}_{n=n_{j-1}+1}P(|X_m|>nj)\geq k.
\end{eqnarray}

Let $h(x)=xg(x),x\geq0$ and $\varphi=h^{-1}$. Then $\varphi$ is a continuous strictly increasing function such that $\varphi(0)=0$ and $\varphi(x)\to\infty$ as $x\to\infty$. Moreover,$\frac{\varphi(x)}{x}=\frac{1}{g(\varphi(x))}$ is strictly decreasing to 0 as $x\to\infty$. By the definitions of the functions $h$ and $\varphi$, and (\ref{3.6}), we have
\begin{eqnarray*}
\lim\limits_{N\to\infty}\inf_{m\geq1}\sum^N_{n=1}P(\varphi(|X_m|)>n)
&=& \lim\limits_{k\to\infty}\inf_{m\geq1}\sum^{n_k}_{n=1}P(\varphi(|X_m|)>n)\\
&=& \lim\limits_{k\to\infty}\inf_{m\geq1}\sum^{n_k}_{n=1}P(|X_m|>h(n))\\
&=& \lim\limits_{k\to\infty}\inf_{m\geq1}\sum^{n_k}_{n=1}P(|X_m|>ng(n))\\
&\geq& \lim\limits_{k\to\infty}k=\infty,
\end{eqnarray*}
which together with  Theorem \ref{thm3.2} and (\ref{eq: equivalent def of WUNI}) implies  that $\{\varphi(|X_m|),m\geq 1\}$ is W-UNI.
\hfill\fbox

Hu and Rosalsky (2015) proved that every nonintegrable nonzero random variable always has slightly more nonintegrability available than that which is assumed. Specifically, they proved the following result which is  an immediate corollary of both  \cite[Theorem 3.3]{CHR16} and the above Theorem \ref{thm3.3}.

\begin{cor} (Hu and Rosalsky (2015, Theorem 1.1)). If $X$ is a random variable with $X\neq 0$ a.s. and $E|X|=\infty$, then there exists a continuous function $g$ on $(0,\infty)$ with $0<g(x)\uparrow \infty$ and $x^{-1}g(x)\downarrow 0$ as $0<x\uparrow\infty$ such that $E(g(|X|))=\infty$.

\end{cor}


\section{Some remarks}

{\bf Remark 4.1.}
In \cite[Remark 5.3]{CHR16}, the authors raised an open problem: whether a UNI sequence of random variables can converge in distribution to an integrable random variable? We can  give a negative answer to this problem. In fact, we can show that {\it any W-UNI sequence of random variables cann't converge in distribution to an integrable random variable.}

Let $\{X_n,n\geq 1\}$ be a W-UNI sequence of random variables, and $\{X_n,n\geq 1\}$ converges in distribution to a random variable $X$. Then by the Skorohod representation theorem, there exist a sequence $\{Y,Y_n,n\geq 1\}$ of random variables such that for any $n\geq 1$, $X_n$ and $Y_n$ has the same distribution,  $X$ and $Y$ have the same distribution, and $\{Y_n,n\geq 1\}$ converges to $Y$ a.s. By the definition of W-UNI, we know that $\{Y_n,n\geq 1\}$ is W-UNI.

For any positive number $a$, by the dominated convergence theorem, we get
$$
E(|Y|\wedge a)=\lim_{n\to\infty}E(|Y_n|\wedge a),
$$
which together with the monotone convergence theorem and the definition of W-UNI, implies that
\begin{eqnarray*}
E(|Y|)&=&\lim_{a\to\infty}E(|Y|\wedge a)\\
&=&\lim_{a\to\infty}\lim_{n\to\infty}E(|Y_n|\wedge a)\\
&\geq&\lim_{a\to\infty}\inf_{n\geq 1}E(|Y_n|\wedge a)=\infty,
\end{eqnarray*}
i.e. $Y$ is NI and thus $X$ is NI.

\bigskip

\noindent{\bf Remark 4.2.} Let $\{X_n,n\geq 1\}$ be a W-UNI sequence of random varialbes. Then

(i) By Definition \ref{defi-1.2}, we know that $\{|X_n|,n\geq 1\}$ is W-UNI,  any subsequence of $\{X_n,n\geq 1\}$ is W-UNI,  and any subsequence of $\{|X_n|,n\geq 1\}$ is W-UNI;

(ii) By (i) and  Remark 4.1, we know that no subsequence of $\{X_n,n\geq 1\}$ can converges  in distribution to a real number;

(iii) By (i) and  Remark 4.1, we know that no subsequence of $\{|X_n|,n\geq 1\}$ can converges  in distribution to a real number;

(iv) By (ii), (iii) and the relation that UNI $\Rightarrow$ W-UNI, we get the first part of \cite[Remark 5.1]{CHR16}.

\bigskip

\noindent {\bf Remark 4.3.} In \cite[Remark 5.2]{CHR16}, the authors showed  that if $\{X_n,n\geq 1\}$ is a UNI sequence of random variables, then no subsequence of $\{|X_n|,n\geq 1\}$ can approach $\infty$ in probability and, a fortiori, no subsequence of $\{X_n,n\geq 1\}$ can approach $\infty$ in probability. Further,   an example is given to show that a sequence of NI random variables can approach $\infty$ in probability. That example is as follows:

Let $Y$ be a NI random variable with $Y>0$ a.s. and let $Y_n=nY,n\geq 1$. Then $\{Y_n,n\geq 1\}$ is a sequence of NI random variable with $Y_n\stackrel{P}{\rightarrow}\infty$.

By Remark \ref{rem1.1}, we know that $\{Y_n,n\geq 1\}$ is W-UNI, and thus W-UNI sequence of random variables can approach $\infty$ in probability.

%
%
%

\bigskip

\noindent{\bf Acknowledgements.}  We thank the support of NNSFC (Grant No. 11371191).

\end{document}